\documentclass{ifacconf}

\usepackage{graphicx}      
\usepackage{natbib}        
\usepackage{gensymb}
\usepackage{epstopdf}
\usepackage{caption}
\usepackage{subcaption}
\usepackage[export]{adjustbox}
\usepackage{amsmath}
\usepackage{xcolor}
\usepackage{soul}
\begin{document}
\begin{frontmatter}

\title{Diagnosing and Decoupling the Degradation Mechanisms in Lithium Ion Cells: An Estimation Approach} 

\author[OSU]{Raja Abhishek Appana} 
\author[OSU]{Faissal El Idrissi} 
\author[OSU]{Prashanth Ramesh}
\author[OSU]{Marcello Canova}
\author[HMC]{Chun Yong Kang}
\author[HMC]{Kimoon Um}

\address[OSU]{Center for Automotive Research, The Ohio State University, Columbus, OH 43210 USA (e-mail: appana.3@osu.edu; elidrissi.2@osu.edu; ramesh.47@osu.edu; canova.1@osu.edu)}

\address[HMC]{Hyundai Motor Company, 150, Hyundaiyeonguso-ro, Namyang-eup, Hwaseong-si, Gyuonggi-do, South Korea (e-mail: chunyongkang@hyundai.com; kimoon.um@hyundai.com)}

\begin{abstract}                
Understanding battery degradation in electric vehicles (EVs) under real-world conditions remains a critical yet under-explored area of research. Central to this investigation is the challenge of estimating the specific degradation modes in aged cells with no available information on usage history, bypassing the conventional yet invasive method of tear-down tests. Using an electrochemical model, this study pioneers a methodology to decouple and isolate the aging mechanisms in batteries sourced from EVs with varying mileages. A robust correlation is established between the model parameters and distinct degradation processes, enabling the diagnosis and estimation of each mechanism's impact on the battery's parameters. This paper sheds light on battery degradation in real-world scenarios and demonstrates the feasibility of their identification, isolation, and approximate quantification of their effects.
\end{abstract}


\end{frontmatter}

\section{Introduction}

The adoption of electric vehicles (EVs) has steadily increased over recent years, underpinned by advancements in energy and power density, safety protocols, longevity, and reduced environmental impact  [\cite{bibra2022global}]. As EVs operate in real-world conditions, their batteries show a gradual decline in capacity, resulting from various degradation phenomena at play. However, estimating the extent of each degradation mechanism in the used battery is challenging, particularly without resorting to destructive testing methods.

Battery aging is a well-documented occurrence, both in controlled laboratory settings and in real-world conditions. \cite{birkl2017degradation} presents experimental evidence of various degradation modes that exist in batteries, spanning the main components (current collectors, electrodes, electrolyte, and separator). Furthermore, research by \cite{zhu2020investigation}, \cite{dubarry2011evaluation}, and \cite{salyer2021extended} consolidates the degradation mechanisms into 3 distinct modes - loss of lithium inventory (LLI), loss of active material (LAM), and impedance increase. LAM is further categorized into loss of active material in the anode ($\textnormal{LAM}_\textnormal{NE}$), and loss of active material in the cathode ($\textnormal{LAM}_\textnormal{PE}$). These mechanisms, though categorized into distinct modes, exhibit a high degree of coupling.

There is extensive research that delves into each mechanism individually, enabling the development of non-destructive diagnostic methods. \cite{dubarry2012synthesize} pioneered a diagnostic approach by artificially simulating the degradation modes and observing the corresponding changes in differential capacity analysis (DCA) and differential voltage analysis (DVA) curves. This method is designed to synthetically generate the cell curves under assumed degradation mechanisms. \cite{birkl2017degradation} developed a diagnostic algorithm utilizing a parameterized open-circuit voltage (OCV) model to determine the degradation modes. By analyzing the changes in open-circuit potentials (OCPs) of the electrodes and mapping pseudo-OCV measurements from the model, they identify degradation mechanisms. \cite{li2019degradation} conducted differential voltage analysis on full cells to identify the aged electrode and estimate capacity loss. \cite{ando2018degradation} devised an empirical curve-fitting methodology to establish the relation between usage conditions and degradation factors.

Some studies proposed the use of electrochemical models to diagnose aging mechanisms and to identify model parameters associated with these phenomena. \cite{uddin2016characterising} developed a technique that identifies the P2D model parameters affecting the calendar aging of a NCA/Graphite cell. \cite{kim2021effective} compared the P2D model calibration parameters of cells at beginning of life (BOL) and end of life (EOL) to identify the parameters influencing aging. \cite{fan2023nondestructive} evaluated cells at different degradation stages and plotted evolution trajectories of aging-related parameters. The electrode operating lithiation ranges and the volume fractions are found to affect the degradation rates. In all these studies, the model is calibrated only on the basis of terminal voltage.

Much of the existing literature on non-destructive battery aging diagnostics relies on data generated in controlled laboratory environments, accompanied by detailed cell cycling history. Only a handful of studies have attempted to decompose the coupled degradation phenomena, using either the DCA or the comparison of electrochemical model calibration parameters. This research focuses on diagnosing batteries extracted directly from actively used EVs. Leveraging physics-based modeling, a diagnostic process is proposed to distinguish and decouple the degradation mechanisms in aged cells, and to quantify their influence on the electrochemical parameters. The study also underscores the necessity of combining IC analysis and electrochemical model calibration to accurately identify these mechanisms. This paper endeavors to bridge the gap between controlled laboratory findings and real-world battery usage, providing a methodological foundation for estimating the health of EV batteries.

\section{Experimental data}

\begin{figure}[hbt!]
    \begin{center}
    \includegraphics[width=8.4cm]{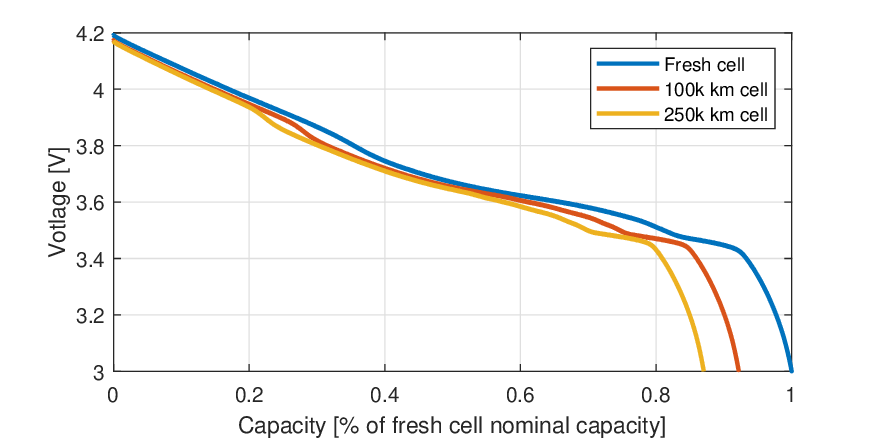}
    \caption{Terminal voltage plots of fresh and aged cells.} 
    \label{fig:OCVs}
    \end{center}
\end{figure}

\begin{figure}[hbt!]
    \begin{center}
    \includegraphics[width=8.4cm]{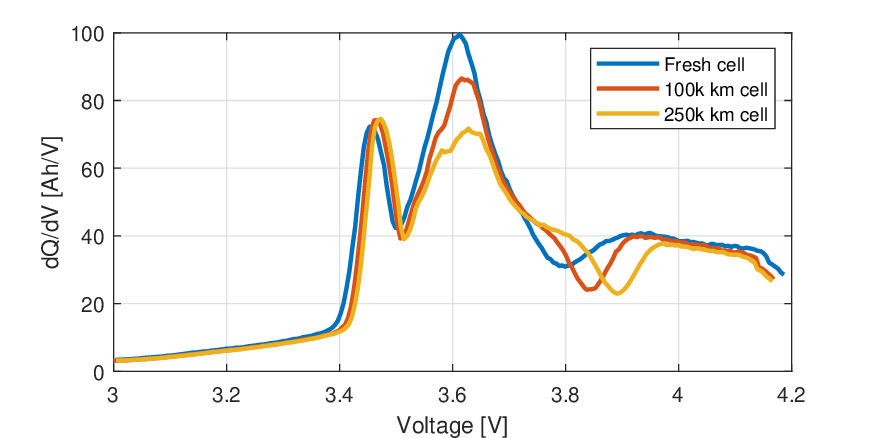}
    \caption{DCA plots of fresh and aged cells.} 
    \label{fig:dQdV}
    \end{center}
\end{figure}

For this study, cells are taken from EVs operating in the real world, with approximate odometer readings of 100,000 km and 250,000 km. The cells were acquired from vehicles after 4-5 years of usage, indicating a significant impact of calendar time and driving throughput on degradation. A fresh cell is also included in testing to establish a baseline. These cells are tested at $25\degree C$ through a capacity test with a C-rate of C/20. Fig. \ref{fig:OCVs} displays the terminal voltage plots of the three cells, with voltage plotted against capacity as a fraction of the nominal capacity of the fresh cell. The cells with mileages of 100k km and 250k km exhibit capacity degradation of $7.80 \%$ and $12.96 \%$ respectively. The fresh cell is also subjected to dynamic tests consisting of charge and discharge pulses to aid in the calibration of the electrochemical model. Fig. \ref{fig:dQdV} illustrates the DCA performed on the C/20 discharge profile of the cells.

\section{Electrochemical model formulation}

\subsection{Extended Single Particle Model (ESPM)}

This study employs the Extended Single Particle Model (ESPM) developed by \cite{fan2016modeling} as the electrochemical model. MATLAB and Simulink are utilized to build the model and run simulations. For conciseness, the equations underlying this model are not presented here. Readers are encouraged to refer to the cited paper for a comprehensive understanding of ESPM. The model was calibrated using experimental data obtained from the fresh cell (NMC - Graphite chemistry) at BOL. The calibrated parameters relevant to this study are given in table \ref{tab:0k}.

\begin{table}[hbt!]
\begin{center}
\caption{Calibrated parameters for a fresh cell.}\label{tab:0k}
\begin{tabular}{ccc}
Parameter & Description & Value \\\hline
$c_{n_0}$  [$mol/m^3$] & Initial lithium ion concentration in NE & 2.75e4 \\ \hline
$c_{p_0}$  [$mol/m^3$] & Initial lithium ion concentration in PE & 2.56e4 \\ \hline
$\varepsilon_n $  [ ]& Active material volume fraction in NE & 0.582 \\ \hline
$\varepsilon_p$  [ ] & Active material volume fraction in PE & 0.540 \\ \hline
\end{tabular}
\end{center}
\end{table}

\subsection{Electrochemical model parameters affecting the aging mechanisms}

\begin{figure*}
    \centering
    \begin{subfigure}[hbt!]{0.49\textwidth}
        \centering
        \includegraphics[width=\textwidth]{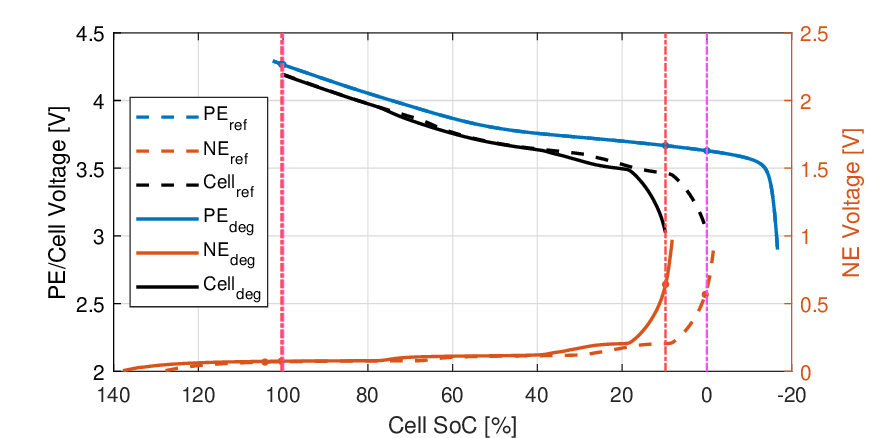}
        \caption{Electrode OCVs with LLI.} 
        \label{fig:10LLIOCVs}
    \end{subfigure}
    \hfill
    \begin{subfigure}[hbt!]{0.49\textwidth}
        \centering
        \includegraphics[width=\textwidth]{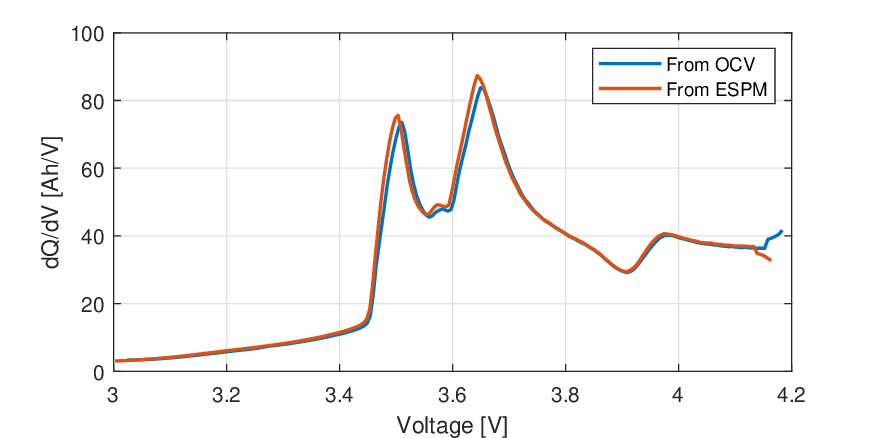}
        \caption{DCA plots with LLI.} 
        \label{fig:10LLIdQdV}
    \end{subfigure}
    \vfill
    \begin{subfigure}[hbt!]{0.49\textwidth}
        \centering
        \includegraphics[width=\textwidth]{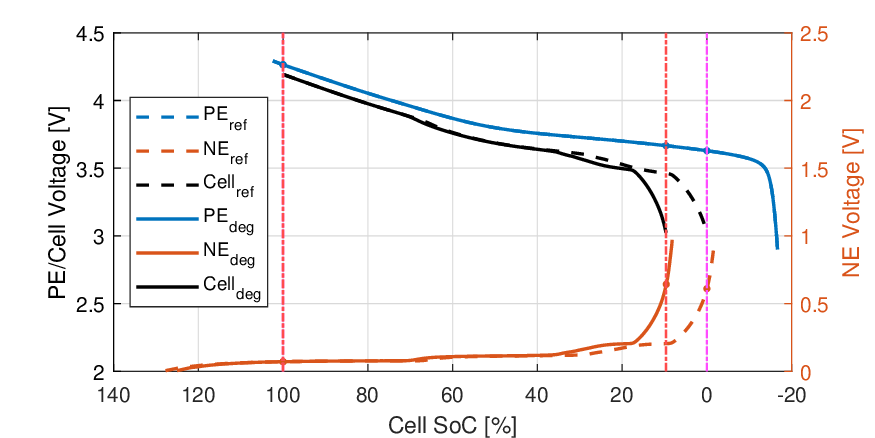}
        \caption{Electrode OCVs with $\textnormal{LAM}_\textnormal{NE}$.} 
        \label{fig:10LAMNEOCVs}
    \end{subfigure}
    \hfill
    \begin{subfigure}[hbt!]{0.49\textwidth}
        \centering
        \includegraphics[width=\textwidth]{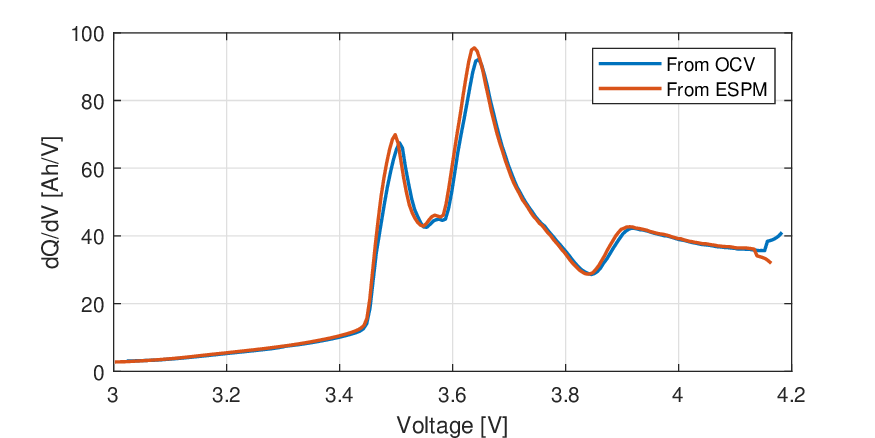}
        \caption{DCA plots with $\textnormal{LAM}_\textnormal{NE}$.} 
        \label{fig:10LAMNEdQdV}
    \end{subfigure}
    \vfill
    \begin{subfigure}[hbt!]{0.49\textwidth}
        \centering
        \includegraphics[width=\textwidth]{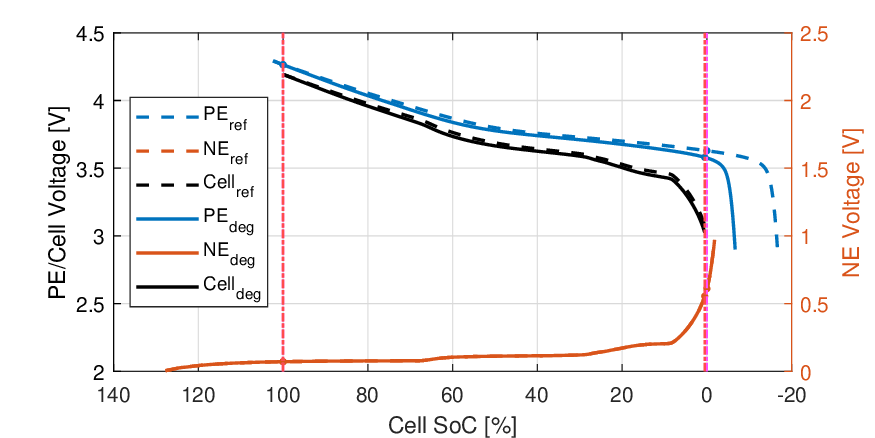}
        \caption{Electrode OCVs with $\textnormal{LAM}_\textnormal{PE}$.} 
        \label{fig:10LAMPEOCVs}
    \end{subfigure}
    \hfill
    \begin{subfigure}[hbt!]{0.49\textwidth}
        \centering
        \includegraphics[width=\textwidth]{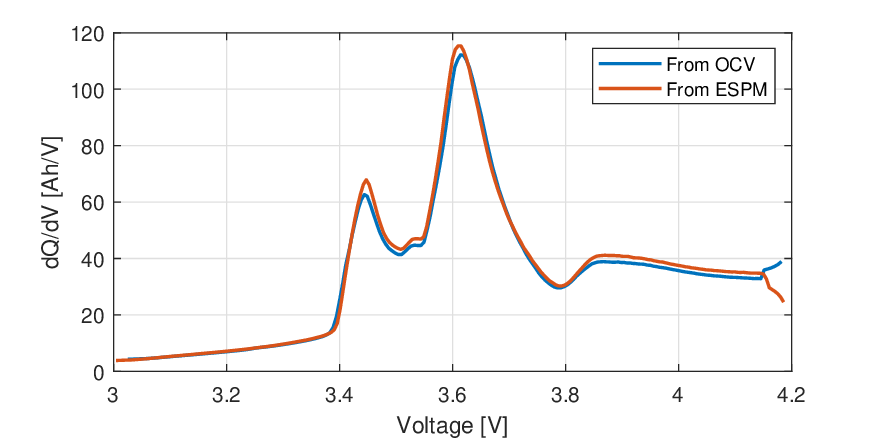}
        \caption{DCA plots with $\textnormal{LAM}_\textnormal{PE}$.} 
        \label{fig:10LAMPEdQdV}
    \end{subfigure}
    \caption{\ref{fig:10LLIOCVs}, \ref{fig:10LAMNEOCVs}, and \ref{fig:10LAMPEOCVs}: Electrode OCP plots. '\textit{ref}' and '\textit{deg}' refer to the fresh state and degraded state curves. \\
    \ref{fig:10LLIdQdV}, \ref{fig:10LAMNEdQdV}, and \ref{fig:10LAMPEdQdV}: DCA plots for each mechanism. '\textit{From OCV}' curve is plotted using the OCV data of aged cells obtained in Fig. \ref{fig:10LLIOCVs}, \ref{fig:10LAMNEOCVs}, and \ref{fig:10LAMPEOCVs}. '\textit{From ESPM}' curve is plotted using the terminal voltage obtained using ESPM after modifying the electrochemical parameters}
    \label{fig:10dQdV}
\end{figure*}

The theoretical capacity of each electrode of the cell $Q_k$ can be calculated using (\ref{eq:Qelec}), where $k$ represents the electrode, $\theta_{0,k}$ and $\theta_{f,k}$ are the initial and final lithium compositions in the electrode, $\varepsilon_k$ is the active material volume fraction, $V_k$ is the volume of active material, and $c_{s,max,k}$ is the lithium saturation concentration.
\begin{equation}\label{eq:Qelec}
    Q_{k} = (\theta_{0,k}-\theta_{f,k})*\varepsilon_k * V_k * c_{s,max,k} * F \quad k =n,p
\end{equation}

\cite{birkl2017degradation} elaborate how the operating stoichiometric regions of the electrodes and their voltage curves change due to the degradation processes. These changes are briefly summarized for convenience. LLI causes a loss of usable Li $x_{LLI}$, due to the occurrence of electrochemical side reactions in the NE. These introduce a side reaction current that leads to the loss of lithium in the electrode, resulting in a relative movement of the NE OCP curve towards the higher potential region of the PE. Consequently, the cathode and anode potentials at the end of discharge (EoD) increase (see Fig. \ref{fig:10LLIOCVs}). This results in a reduction of the normalized initial lithium concentration in NE $\theta_{0,n}$.
\begin{equation}\label{eq:QNELLI}
    Q_{k,aged} = (\theta_{0,k}-\theta_{f,k}-x_{LLI})*\hat{\varepsilon}_k * V_k * c_{s,max,k} * F \quad k =n,p
\end{equation}

The initial lithium ion concentration in the negative electrode $c_{n_0}$ can be calculated from the initial lithium composition by using (\ref{eq:cn0}), showing that changes due to NE side reactions can be captured by $c_{n_0}$. This is also verified via the electrochemical model. $c_{n_0}$ is calibrated to achieve the same capacity degradation obtained in Fig. \ref{fig:10LLIOCVs} due to LLI. The DCA curves that are plotted using the theoretically calculated OCV (from Fig. \ref{fig:10LLIOCVs}) and the terminal voltage obtained from the ESPM are compared in Fig. \ref{fig:10LLIdQdV}, showing a good match of the results.
\begin{equation} \label{eq:cn0}
c_{n_0} = c_{max_n} * \theta_{0,n}
\end{equation}

From the perspective of electrochemical modeling, it is generally assumed that the effect of LAM is evident when NE is lithiated and PE is delithiated (\cite{dubarry2022best}). These changes are illustrated in Fig. \ref{fig:10LAMNEOCVs}-\ref{fig:10LAMPEOCVs}. The OCV curve of the electrode with LAM shrinks around the lithium composition at $100\%$ cell SoC, which indirectly reduces $\theta_f$ in the affected electrode. Comparing with (\ref{eq:QNELLI}), $x_{LLI}$ accounts for the shrinking of the lithiation windows and $\hat{\varepsilon}_k$ reduces from the original value due to a loss in accessible lithium.

\cite{carnovale2020modeling} model LAM in each electrode using a relation between the rate of electrode active material volume fraction $\varepsilon$ and the electrode current. So, the negative electrode volume fraction $\varepsilon_n$ and positive electrode volume fraction $\varepsilon_p$ can be directly attributed to $\textnormal{LAM}_\textnormal{NE}$ and $\textnormal{LAM}_\textnormal{PE}$. These parameters are reduced individually in the electrochemical model to simulate the capacity loss obtained in Fig. \ref{fig:10LAMNEOCVs} and \ref{fig:10LAMPEOCVs} respectively. Reducing these parameters directly affects the intercalation current density in the solid phase, which in turn affects the accessible lithiation ranges in the electrochemical model, as the cell is charged and discharged between fixed terminal voltage limits. The DCA curves plotted using the output voltage data of the ESPM are compared with the DCA curves obtained earlier using OCVs. The results are in agreement, verifying the hypothesis made (see Fig.\ref{fig:10LAMNEdQdV}-\ref{fig:10LAMPEdQdV}). 

The analysis conducted in this section confirms that estimating the parameters $c_{n_0}, \varepsilon_n, \varepsilon_p$ from data of aged cells can be used to isolate and quantify the degradation mechanisms in the positive and negative electrode and predict their effects on cell voltage using any electrochemical model.

\section{Results and Discussion}

This section outlines the methodology developed to diagnose the degradation mechanisms and estimate the electrochemical parameters leveraging the ESPM. The quantitative effect of each degradation mechanism on cell characteristics is highly dependent on the cell history (applied current, calendar life, etc.). Firstly, the necessity of a degradation mechanism to explain the behavior of the aged cell is established. Secondly, if a degradation mechanism is indeed necessary, the sufficiency of that mechanism alone to explain aged cell behavior is verified. Once these questions are addressed and degradation mechanisms are identified, the parameters of the ESPM for the aged cell are estimated to quantify the effect of degradation mechanisms and reproduce the voltage characteristics of the aged cell.

\begin{figure}[hbt!]
    \begin{center}
    \includegraphics[width=8.4cm]{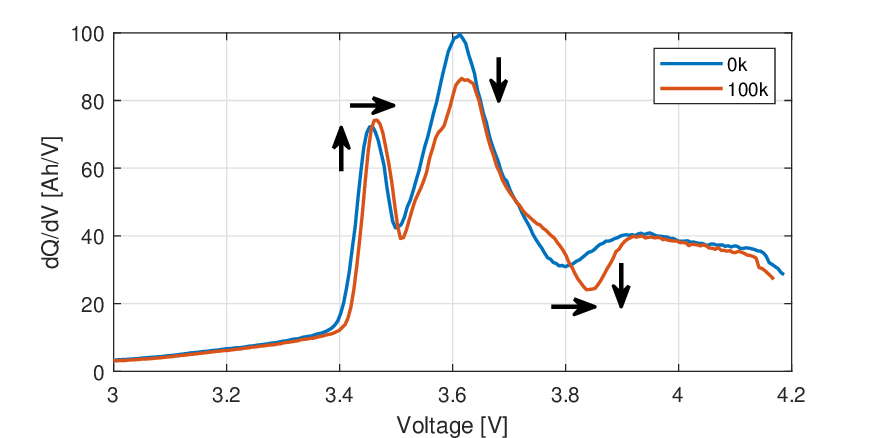}
    \caption{DCA plot of fresh cell and 100k km cell using experimentally measured data.}
    \label{fig:dQdVexp}
    \end{center}
\end{figure}

Fig. \ref{fig:dQdVexp} compares the DCA plots of a fresh cell and an aged cell harvested from the battery pack of an EV driven for 100k km. Arrows indicate the relative movement of features in the DCA curve. For instance, in the aged cell, the low voltage peak at $3.45V$ shifts slightly towards higher voltage and intensifies in magnitude, while the high voltage valley at $3.8V$ shifts significantly towards a higher potential with a change in magnitude. These indicators are used to identify and decouple the degradation mechanisms. To understand the effect of each mechanism on the DCA curve, the corresponding model parameter is varied while the others remain at their initial values.

\subsubsection{Loss of Lithium Inventory (LLI) - Estimating $c_{n_0}$:}
\begin{figure}[hbt!]
    \begin{center}
    \includegraphics[width=8.4cm]{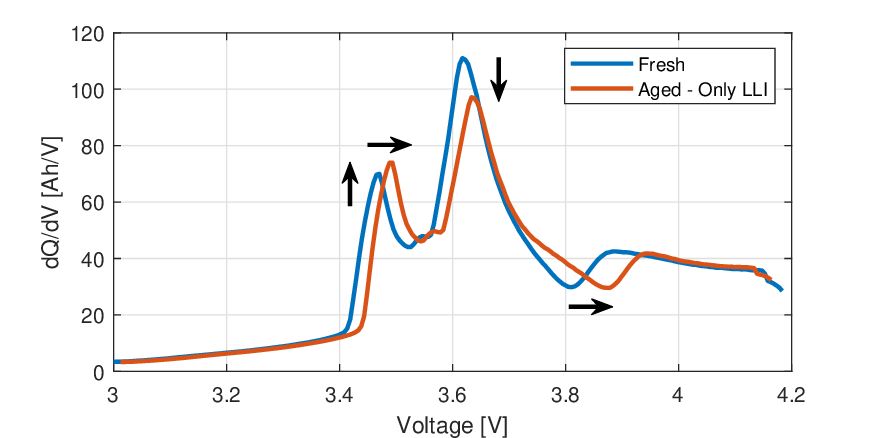}
    \caption{DCA plot comparing the simulated output of fresh cell model and aged cell model with LLI.}
    \label{fig:dQdVLLI}
    \end{center}
\end{figure}

As established earlier, varying the initial lithium-ion concentration in the negative electrode $c_{n_0}$ within the electrochemical model simulates LLI. Initially, $c_{n_0}$ is calibrated in the ESPM to accurately match the simulated capacity with the measured 100k km cell capacity. The DCA curve from this simulation does not show good coherence with the experimental data. The high voltage valley (at $3.8V$) does not exhibit a decrease in intensity similar to the real-world sample, showing that LLI alone is insufficient to explain the aged cell behavior. Furthermore, based on Fig. \ref{fig:dQdVLLI}, it is evident that only LLI increases the low voltage peak's (at $3.45 V$) intensity and induces the pronounced shift in the high voltage valley (at $3.8 V$) to the right. These specific characteristics cannot be replicated with $\textnormal{LAM}_\textnormal{NE}$ and $\textnormal{LAM}_\textnormal{PE}$ mechanisms (refer to Fig. \ref{fig:dQdVLAMNE}, and \ref{fig:dQdVLAMPE}). Therefore, the definite presence of LLI in the aged cell can be established.

\begin{figure}[hbt!]
    \begin{center}
    \includegraphics[width=8.4cm]{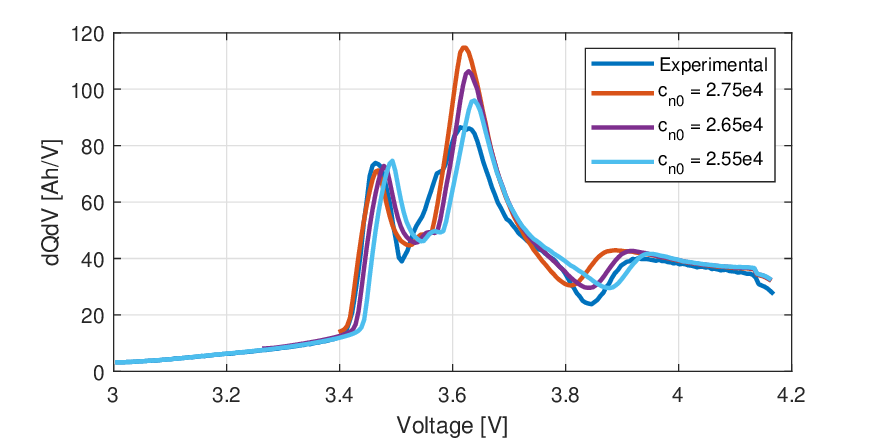}
    \caption{DCA plots with various $c_{n_0}$ - Experimental vs Simulated - 100k cell.}
    \label{fig:dQdV_Cn0}
    \end{center}
\end{figure}

By observing Fig. \ref{fig:dQdVLLI}, \ref{fig:dQdVLAMNE}, and \ref{fig:dQdVLAMPE}, it can also be said that the valley at $3.8V$ shows the most distinguishable characteristic i.e. it shifts significantly in the presence of LLI. This feature is leveraged to estimate $c_{n_0}$, the parameter affecting LLI. The objective of this step is to match the position of the high-voltage valley, using the loss function in (\ref{eq:loss1}). $V_{deg,HV valley}$ and $V_{ref,HV valley}$ are the positions of the high-voltage valley in the aged cell and fresh cell simulations respectively. A simple grid search is performed with various values of $c_{n_0}$ (see Fig. \ref{fig:dQdV_Cn0}) to minimize the objective function. All the grid search points are not shown in the figure. While the final estimated value, $c_{n_0} = 2.65e4 mol/m^3$, matches the position of the high voltage valley, the model predicted capacity is higher than the measured capacity of the aged cell indicating the need to include additional degradation mechanisms. 
\begin{equation} \label{eq:loss1}
    J = |V_{deg,HV valley} - V_{ref,HV valley}|
\end{equation}

\subsubsection{Loss of Active material in negative electrode ($\textnormal{LAM}_\textnormal{NE}$) - Estimating $\varepsilon_n$:}
\begin{figure}[hbt!]
    \begin{center}
    \includegraphics[width=8.4cm]{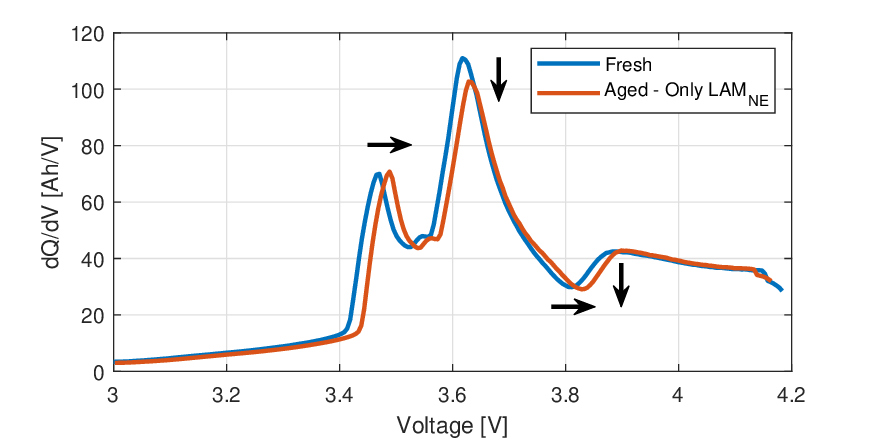}
    \caption{DCA plot comparing the simulated output of fresh cell model and aged cell model with $\textnormal{LAM}_\textnormal{NE}$.}
    \label{fig:dQdVLAMNE}
    \end{center}
\end{figure}

At this stage, it is established that LLI is present in the aged cell. But, to verify the necessity of $\textnormal{LAM}_\textnormal{NE}$, only the NE active material volume fraction $\varepsilon_n$ is adjusted to match the simulated terminal capacity with the aged cell capacity, while keeping $c_{n_0}$ at the BOL calibration value. Comparison of the simulated DCA curves of aged and fresh cells (see Fig. \ref{fig:dQdVLAMNE}) reveals that the decrease in intensity of the high voltage valley (at $3.8V$) is exclusively attributable to $\textnormal{LAM}_\textnormal{NE}$, not to LLI or $\textnormal{LAM}_\textnormal{PE}$ (refer to Fig. \ref{fig:dQdVLLI}, and \ref{fig:dQdVLAMPE}). Consequently, $\textnormal{LAM}_\textnormal{NE}$ is needed to explain cell behavior. However, this mechanism does not amplify the intensity of the low voltage peak (at $3.45V$), which is a characteristic of the aged cell as seen in Fig. \ref{fig:dQdVexp}. Therefore, it points to the presence of additional mechanisms beyond $\textnormal{LAM}_\textnormal{NE}$ to account for this behavior.

\begin{figure}[hbt!]
    \begin{center}
    \includegraphics[width=8.4cm]{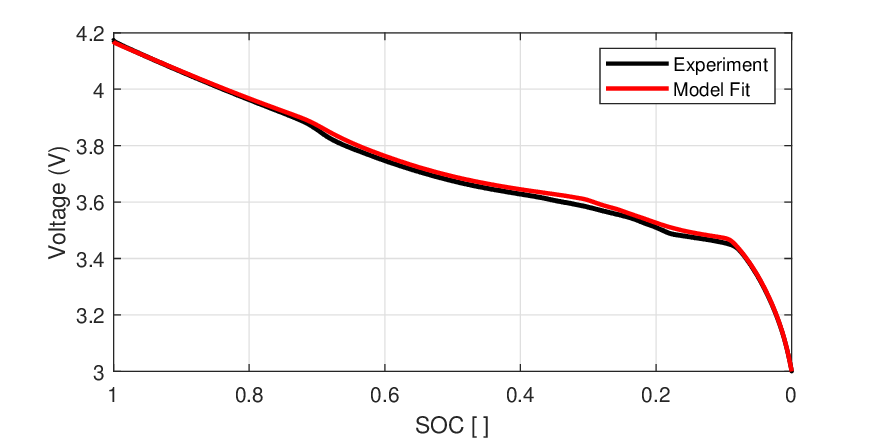}
    \caption{Terminal voltage plot of 100k km aged cell - Experimental vs Simulated with LLI + $\textnormal{LAM}_\textnormal{NE}$.}
    \label{fig:Voltage100kLLILAMNE}
    \end{center}
\end{figure}

\begin{figure}[hbt!]
    \begin{center}
    \includegraphics[width=8.4cm]{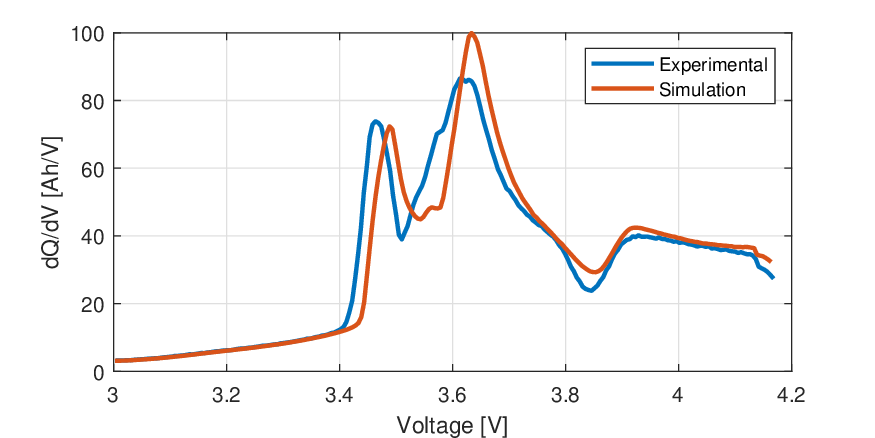}
    \caption{DCA plot of 100k km aged cell - Experimental vs Simulated with LLI + $\textnormal{LAM}_\textnormal{NE}$.}
    \label{fig:dQdVLLILAMNE}
    \end{center}
\end{figure}

With both LLI and $\textnormal{LAM}_\textnormal{NE}$ confirmed, and with an already estimated value of $c_{n_0}$, $\varepsilon_n$ is tuned to minimize the loss function in (\ref{eq:loss2}). The value of $\varepsilon_n = 0.561$ is chosen to match the estimated terminal capacity with the measured value. The RMS error on the terminal voltage plot (see Fig. \ref{fig:Voltage100kLLILAMNE}) is $14mV$ and the DCA plot is shown in Fig. \ref{fig:dQdVLLILAMNE}. The plot shows a slight mismatch in the peak positions (mainly at $3.8V$) due to the shift induced from $\textnormal{LAM}_\textnormal{NE}$.
\begin{equation} \label{eq:loss2}
    J = \sqrt{(\Sigma(V_{deg}-V_{ref})^2)/n}
\end{equation}

\subsubsection{Loss of Active material in positive electrode ($\textnormal{LAM}_\textnormal{PE}$) - Estimating $\varepsilon_p$:}
\begin{figure}[hbt!]
    \begin{center}
    \includegraphics[width=8.4cm]{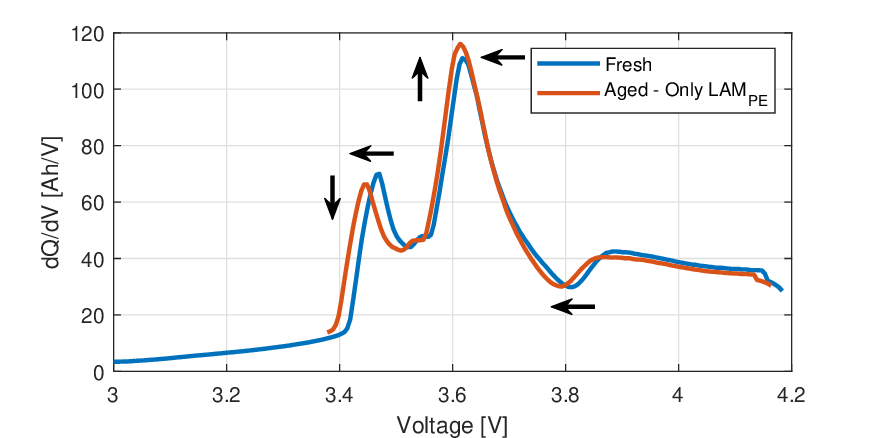}
    \caption{DCA plot comparing the simulated output of fresh cell model and aged cell model with $\textnormal{LAM}_\textnormal{PE}$.}
    \label{fig:dQdVLAMPE}
    \end{center}
\end{figure}

The calibration obtained at the previous step serves as a good starting point to examine whether $\textnormal{LAM}_\textnormal{PE}$ is necessary. In Fig. \ref{fig:dQdVLLILAMNE}, the notable shift of the valley at $3.8V$ towards higher voltages compared to the 100k km aged cells is evident. This necessitates the inclusion of $\textnormal{LAM}_\textnormal{PE}$, as this mechanism counteracts this shift. Since the negative electrode serves as the capacity-limiting electrode for the cell, the presence of $\textnormal{LAM}_\textnormal{PE}$ does not reduce cell capacity (also observed in Fig. \ref{fig:10LAMPEOCVs}). Additionally, $\textnormal{LAM}_\textnormal{PE}$ shifts the DCA curve towards lower voltages (see Fig. \ref{fig:dQdVLAMPE}), which starkly contrasts with the aged cell. This implies that $\textnormal{LAM}_\textnormal{PE}$ cannot be present as the only degradation mechanism in the cell. In summary, all three mechanisms need to be present together to explain the aged cell behavior. 

\begin{figure}[hbt!]
    \begin{center}+
    \includegraphics[width=8.4cm]{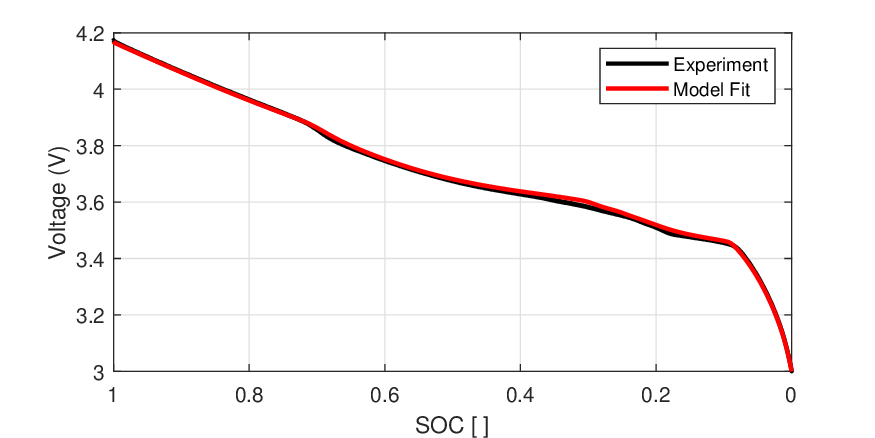}
    \caption{Calibrated model terminal voltage plot of 100k km cell.}
    \label{fig:Voltage100k}
    \end{center}
\end{figure}

\begin{figure}[hbt!]
    \begin{center}
    \includegraphics[width=8.4cm]{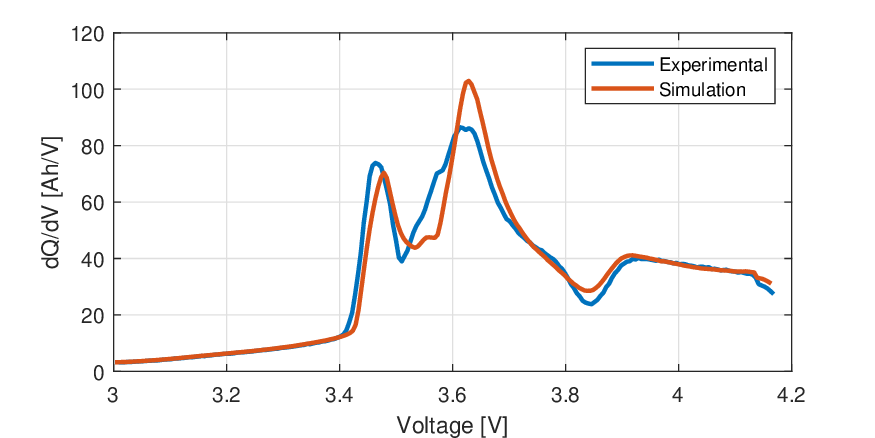}
    \caption{Experimental and Simulated DCA plots of 100k km cell.}
    \label{fig:dQdV100k}
    \end{center}
\end{figure}

Now, $\varepsilon_p$ is calibrated using the same loss function in (\ref{eq:loss2}) to negate that shift and match the peak positions. Results of the calibration process are given in table \ref{tab:100k}. The RMS error on the simulated voltage is $7.96 mV$. The terminal voltage plot and the DCA curve are shown in Fig. \ref{fig:Voltage100k} and \ref{fig:dQdV100k}. Though the quantitative effect of each mechanism on capacity cannot be exactly determined without destructive tests, their effect on the cell characteristics can be determined by this method. It is seen that all the parameters have decreased by around $3.6\%$.

\begin{table}[hbt!]
\begin{center}
\caption{Aged cell - Calibrated electrochemical parameters at 100k km.}\label{tab:100k}
\begin{tabular}{ccc}
\textbf{Parameter} & \textbf{BOL value}  & \textbf{Value at 100k}  \\\hline
$c_{n_0}$ [$mol/m^3$]& $2.75e4$& $2.65e4$\\ \hline
$\varepsilon_n$ [ ]& $0.582$ & $0.561$ \\ \hline
$\varepsilon_p$ [ ]& $0.540$ & $0.520$ \\ \hline
\end{tabular}
\end{center}
\end{table}

\begin{figure}[hbt!]
    \begin{center}
    \includegraphics[width=8.4cm]{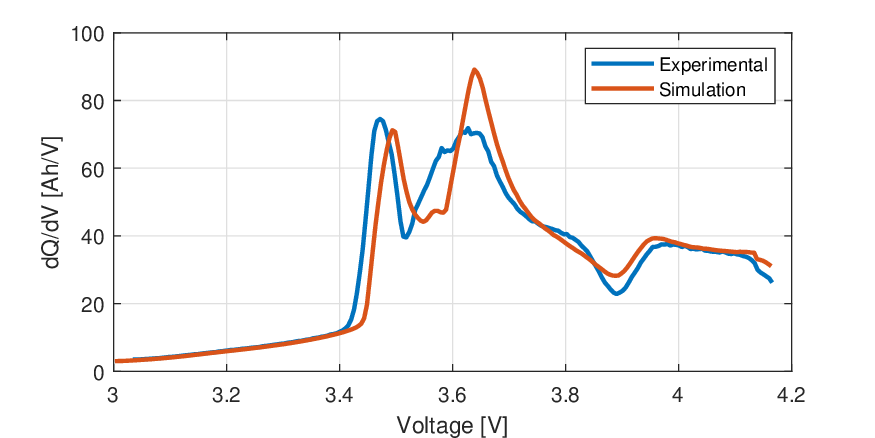}
    \caption{Experimental and Simulated DCA plots of 250k km aged cell.}
    \label{fig:dQdV250k}
    \end{center}
\end{figure}

This methodology is also applied to the 250k km cell to identify the aging mechanisms and to estimate the electrochemical parameters for an additional cell sample. All three aging mechanisms are identified based on the experimental data and the corresponding electrochemical parameters are estimated (see table \ref{tab:250k}). The DCA plot from the resulting calibration is shown in Fig. \ref{fig:dQdV250k}.

\begin{table}[hbt!]
\begin{center}
\caption{Aged cell - Calibrated electrochemical parameters at 100k and 250k km.}\label{tab:250k}
\begin{tabular}{cccc}
\textbf{Parameter} & \textbf{BOL value} & \textbf{100k} & \textbf{250k} \\\hline
$c_{n_0}$ [$mol/m^3$]& $2.75e4$& $2.65e4$ & $2.515e4$\\ \hline
$\varepsilon_n$ [ ]& $0.582$ & $0.561$ & $0.557$ \\ \hline
$\varepsilon_p$ [ ]& $0.540$ & $0.520$ & $0.515$ \\ \hline
\end{tabular}
\end{center}
\end{table}

\section{Conclusions and Future work}

Through this study, we have established a correlation between various degradation mechanisms and electrochemical parameters using fundamental principles. Specifically, LLI, $\textnormal{LAM}_\textnormal{NE}$, and $\textnormal{LAM}_\textnormal{PE}$ directly influence the ESPM parameters $c_{n_0}, \varepsilon_n, \textnormal{and}, \varepsilon_p$ respectively. This correlation is leveraged to identify and isolate aging mechanisms in cells extracted from EVs used in the real world with mileages of 100k and 250k km. The developed methodology is capable of diagnosing different degradation mechanisms through observations based on DCA characteristics. Furthermore, the study highlights the necessity of integrating DCA with electrochemical modeling to effectively identify these mechanisms. A simple estimation technique is developed and applied to estimate various electrochemical parameters in an aged cell. The results can be further used to evaluate the contribution of each degradation mechanism towards capacity loss. Finally, our findings conclusively demonstrate the presence of all three considered mechanisms in real-world EV batteries. The team is currently working on validating the results using post-mortem characterizations and aging model simulations.

\begin{ack}
The OSU team acknowledges and thanks Hyundai Motor Company for funding support.
\end{ack}

\bibliography{ifacconf}             

\end{document}